\documentclass[11pt,reqno]{amsart}
\ifdefined\pdfglyphtounicode
  \input{glyphtounicode}
\fi
\usepackage[T1]{fontenc}
\usepackage[utf8]{inputenc}
\usepackage{lmodern}
\usepackage[a4paper,margin=27mm]{geometry}
\usepackage{amsmath,amssymb,amsthm,mathtools}
\usepackage{microtype}
\usepackage{xcolor}
\definecolor{linkblue}{RGB}{25,60,105}
\usepackage[colorlinks=true,linkcolor=linkblue,citecolor=linkblue,urlcolor=linkblue]{hyperref}
\usepackage{bookmark}
\allowdisplaybreaks[2]
\numberwithin{equation}{section}
\newtheorem{theorem}{Theorem}[section]
\newtheorem{lemma}[theorem]{Lemma}
\newtheorem{proposition}[theorem]{Proposition}
\newtheorem{corollary}[theorem]{Corollary}
\title[Frustration of hypercubes]{The frustration index of fully frustrated hypercubes}
\author{Lujia Wang}
\address{Hainan Bielefeld University of Applied Sciences, Danzhou, China}
\email{lujia.wang@hibiuh.edu.cn}
\date{11 November 2026}
\hypersetup{pdftitle={The frustration index of fully frustrated hypercubes},pdfauthor={Lujia Wang}}
\begin{document}
\begin{abstract}
The frustration index of a signed graph is the minimum number of edges
meeting every negative cycle. Let $F_n$ be the frustration index of the
signed $n$-cube $Q_n$ in which every square is negative. Using signings
constructed from weighted cubes, we prove that for every positive integer $r$,
\[
 F_{r^2+1}=(r^2-r+1)2^{r^2-1}.
\]
When $4\mid r$, our construction also gives the exact value
$F_n=(n-r-p/4)2^{n-2}$ in the six dimensions $n=r^2+1+pr/2$,
$p\in\{\pm1,\pm2,\pm3\}$. In dimension nine an argument with integral
eigenvectors gives $F_9=772$.

By varying one weight in the construction for $r^2+1$, we obtain
$F_n=(n-\sqrt n+O(n^{-1/2}))2^{n-2}$ for arbitrary $n$. The values
at $n=r^2+1$ show that the error term has the best possible order.
Taking the positive edges of an optimal signing gives a square-free
subgraph of $Q_n$ with at least $(n+\sqrt n-O(n^{-1/2}))2^{n-2}$ edges
whose complement is also square-free.
\end{abstract}
\maketitle

\section{Introduction}
Erd\H{o}s~\cite[Problem~10, p.~8]{Erdos} asked whether every square-free
subgraph of the $n$-cube $Q_n$ has at most $(1/2+o(1))n2^{n-1}$ edges
as $n\to\infty$. One way to construct a large square-free subgraph is
to color the edges of $Q_n$ red and blue so that each square has one
or three red edges. Both color classes are then square-free. If $k$
edges are red, the blue subgraph has $n2^{n-1}-k$ edges. How small
can $k$ be?

Give the red edges sign $-1$ and the blue edges sign $+1$. The condition
on the coloring says precisely that the signing is \emph{fully frustrated},
that is, the product of the edge signs on each square is $-1$.
The \emph{frustration index} of a signed graph is the minimum number of
edges meeting every negative cycle, or equivalently, the minimum number
of negative edges after switching. Here switching at a vertex reverses
the signs of its incident edges. Since the squares generate the cycle
space of $Q_n$, all fully frustrated signings of $Q_n$ differ only by
switching. Their common frustration index $F_n$ is therefore the
minimum possible number of red edges. We call a fully frustrated
signing with $F_n$ negative edges \emph{optimal}.

Bialostocki~\cite{Bialostocki} proved that in any two-coloring of
$E(Q_n)$ with no monochromatic square, each color class has at most
$(n+\sqrt n)2^{n-2}$ edges. Applying this to the blue edges gives
$F_n\ge(n-\sqrt n)2^{n-2}$. Brass, Harborth and Nienborg~\cite{BHN}
constructed colorings attaining Bialostocki's bound when $n$ is a
power of four. Laplante et al.~\cite[Theorem~5.4]{LNSW} showed that
the construction in~\cite{BHN} gives a fully frustrated signing in every
dimension.

Parity can improve this lower bound. The \emph{signed degree} of a
vertex is the number of incident
positive edges minus the number of incident negative edges. Every signed
degree has the parity of $n$. In a fully frustrated signing the mean
square of the signed degrees is $n$, while their sum is $n2^n-4k$ if
$k$ edges are negative. For a positive integer $r$ of parity opposite
to $n$, each signed degree differs from $r$ by at least one. Averaging
the squares of these differences, as in the argument of Derrida
et al.~\cite[pp.~620--621]{DPTV}, gives
\begin{equation}\label{eq:parity-bound}
 F_n\ge\left(n-\frac{n+r^2-1}{2r}\right)2^{n-2}.
\end{equation}
Equality holds precisely when some fully frustrated signing has signed
degree $r-1$ or $r+1$ at every vertex (Lemma~\ref{lem:exact-local-field}).

In dimension $n=r^2+1$, \eqref{eq:parity-bound} reduces to
$F_n\ge(n-r)2^{n-2}$. The construction in~\cite{BHN,LNSW} attains
this bound when $r$ is a power of two or $r=3$. Our main theorem
extends equality to every positive integer $r$.

\begin{theorem}\label{thm:exact-values}
For every positive integer $r$,
\begin{equation}\label{eq:exact-family}
 F_{r^2+1}=(r^2-r+1)2^{r^2-1}.
\end{equation}
If $4\mid r$ and $p\in\{\pm1,\pm2,\pm3\}$, then
\begin{equation}\label{eq:additional-exact}
 F_n=(n-r-p/4)2^{n-2},\qquad n=r^2+1+pr/2.
\end{equation}
In both cases equality holds in \eqref{eq:parity-bound}.
\end{theorem}
For example, $r=4$ gives $F_n$ in every odd dimension from $11$ to $23$.
The equality condition also determines how often the two signed degrees
occur. In any optimal signing of $Q_{r^2+1}$, half the vertices in
each parity class have signed degree $r+1$ and half have signed degree
$r-1$. Here the parity of a vertex in $\{0,1\}^n$ is the parity of its
number of ones. In \eqref{eq:additional-exact} the corresponding
proportions are $\tfrac12+p/8$ and $\tfrac12-p/8$.

We obtain these signings by lifting from weighted cubes while preserving
the distribution of signed degrees (Theorem~\ref{lem:short-unit-lifting}).
On $Q_{r^2+1}$ the two degrees label the vertices by a Boolean
function. For odd $r\ge3$, we can choose the lifting so that this
function changes on exactly half the edges in each direction
(Proposition~\ref{prop:field-directions}). For $r=2^k\pm1\ge3$, we
can also require that every vertex have equally many neighbors of
each signed degree (Corollary~\ref{cor:field-every-vertex}).

To compare $F_n$ across dimensions, let $A_n$ be the signed adjacency
matrix of a fully frustrated signing of $Q_n$. A Boolean function
$f\colon\{0,1\}^n\to\{-1,1\}$ specifies a switching at the vertices
where $f(x)=-1$. The sign of the edge $xy$ is then multiplied by
$f(x)f(y)$, so $f^{\mathsf T}A_nf$ counts each edge twice with its
switched sign. Maximizing over $f$ and using $A_n^2=nI$ gives
\begin{equation}\label{eq:normalization}
 \mu_n=\max_f 2^{-n}f^{\mathsf T}A_nf,\qquad
 F_n=(n-\mu_n)2^{n-2},\qquad \mu_n\le\sqrt n.
\end{equation}

The upper bound of Laplante et al.~\cite[Theorem~5.4]{LNSW} takes
the uniform form
\[
 F_n\le\left(n-\frac{2\sqrt2}{3}\sqrt n+O(n^{-1/2})\right)2^{n-2}.
\]
Along $n=2\cdot4^k$, their lower estimate for $\mu_n$ divided by
$\sqrt n$ tends to $2\sqrt2/3$. Theorem~\ref{thm:exact-values} gives
$\mu_{r^2+1}=r$. Since $\mu_n$ is nondecreasing by
\eqref{eq:monotone}, it follows that $\mu_n\ge\lfloor\sqrt{n-1}\rfloor$
for $n\ge2$, and hence $\sqrt n-\mu_n=O(1)$. By varying one weight
in the same construction, we obtain the following sharper estimate.

\begin{theorem}\label{thm:main}
There is an absolute constant $C$ such that for every $n\ge1$,
\[
 (n-\sqrt n)2^{n-2}\le F_n
 \le\left(n-\sqrt n+\frac{C}{\sqrt n}\right)2^{n-2}.
\]
\end{theorem}
The estimate $0\le\sqrt n-\mu_n\le C/\sqrt n$ has the best possible
order. Indeed, Theorem~\ref{thm:exact-values} gives
$\sqrt{r^2+1}-\mu_{r^2+1}=1/(\sqrt{r^2+1}+r)$.

Write $\operatorname{ex}(Q_n,C_4)$ for the largest number of edges
in a square-free subgraph of $Q_n$. Since every square has edges of
both signs, the $n2^{n-1}-F_n$ positive edges of an optimal signing
form a square-free subgraph whose complement is also square-free.

\begin{corollary}\label{cor:extremal}
Let $C$ be as in Theorem~\ref{thm:main}. For every $n\ge1$ there
is a subgraph $G$ of $Q_n$ such that both
$G$ and its edge complement are $C_4$-free and
\[
 |E(G)|\ge\left(n+\sqrt n-\frac C{\sqrt n}\right)2^{n-2}.
\]
In particular, the right-hand side is a lower bound on
$\operatorname{ex}(Q_n,C_4)$. In the dimensions of
\eqref{eq:exact-family} and \eqref{eq:additional-exact}, respectively,
the construction gives
\[
 \operatorname{ex}(Q_{r^2+1},C_4)\ge (r^2+r+1)2^{r^2-1},
 \qquad
 \operatorname{ex}(Q_n,C_4)\ge (n+r+p/4)2^{n-2}.
\]
\end{corollary}

In square dimensions one can ask when $F_n$ attains the spectral
lower bound $(n-\sqrt n)2^{n-2}$. This question arose in the fully
frustrated Ising model studied by Derrida, Pomeau, Toulouse and
Vannimenus~\cite{DPTV}. With couplings of magnitude $1/\sqrt n$, the
ground-state energy per vertex is $-\mu_n/(2\sqrt n)$. Parisi and
Potters~\cite[Section~1.3]{PP} conjectured that this energy tends to
$-1/2$ and that it equals $-1/2$ whenever $n$ is a square.
Theorem~\ref{thm:main} gives the limit with error $O(n^{-1})$.

The energy equals $-1/2$ exactly when the spectral bound is attained.
This bound is integral precisely when $n$ is a square, and it is
attained precisely when
$n$ is a power of four (Proposition~\ref{prop:square-gap}). The first
square dimension in which it is strict is $9$, where it gives $768$.
Marinari, Parisi and Ritort~\cite[p.~329]{MPR} reported an energy
corresponding to $772$ negative edges, and an explicit switching with
this value is available~\cite{Q9witness}. We prove the matching lower
bound.

\begin{theorem}\label{thm:nine}
The frustration index of the fully frustrated $9$-cube is $F_9=772$.
\end{theorem}
It follows that
the largest edge set of $Q_9$ meeting every square in an odd number
of edges has $9\cdot2^8-772=1532$ edges. This answers a question of
Minamoto~\cite[Section~10, Problem~7]{Minamoto}, who asked whether
$1536$ such edges were possible.

Section~\ref{sec:refinements} develops the lifting used for the exact
families in Section~\ref{sec:exact} and the uniform bound in
Section~\ref{sec:asymptotic}. Section~\ref{sec:q9} proves
$F_9=772$ using integral eigenvectors. We return to square-free
color classes in Section~\ref{sec:conclusion}.

\section{Weighted cubes and refinements}\label{sec:refinements}
Fix a positive integer $r$ and consider a square with both horizontal
edges positive and the two
vertical edges of opposite signs. Give the horizontal edges weight
$r$ and the vertical edges weight $1$. The sum of the incident edge
signs, counted with these weights, is $r-1$ at two vertices and
$r+1$ at the other two. These are just the signed degrees required
for equality in \eqref{eq:parity-bound} on $Q_{r^2+1}$. How can we
realize them on a cube in which all edges have weight one?

For $x=(x_0,\ldots,x_{d-1})$, write $|x|=\sum_jx_j$ for its Hamming
weight. Let $e_i$ be the $i$th coordinate vector and write $\oplus$ for
coordinatewise addition modulo two, so $x\oplus e_i$ is the neighbor
obtained by flipping bit $i$. We use the \emph{canonical signing} in which
the edge from $x$ to $x\oplus e_i$ has sign $(-1)^{\sum_{j<i}x_j}$.
The signed adjacency operator $\Gamma_i$ in direction $i$ transfers the
value at this neighbor to $x$ with the sign of their edge. The parity
operator $P_d$ multiplies $f$ by $1$ on even-weight vertices and by
$-1$ on odd-weight vertices. Thus
\[
 (\Gamma_i f)(x)=(-1)^{\sum_{j<i}x_j}f(x\oplus e_i),\qquad
 (P_df)(x)=(-1)^{|x|}f(x).
\]
The operator $P_d$ records the sign contributed by these $d$ bits to
any coordinate that follows them. For a vector $w=(w_0,\ldots,w_{d-1})$
of real weights, the weighted adjacency operator is
\[
 \Gamma_w=\sum_{i=0}^{d-1} w_i\Gamma_i.
\]
Taking every weight equal to one gives $\Gamma_{(1,\ldots,1)}=A_d$.
After switching by a Boolean function $f$, the weighted signed degree
at $x$ is called its \emph{local field}, denoted by $h(x)$. Its average
over the vertices is the \emph{mean local field}. With the normalized
inner product $\langle f,g\rangle=2^{-d}\sum_x f(x)g(x)$ and the
abbreviation $\langle B\rangle_f=\langle f,Bf\rangle$, we have
\[
 h(x)=f(x)(\Gamma_w f)(x),\qquad
 \langle\Gamma_w\rangle_f=\langle f,\Gamma_w f\rangle=2^{-d}\sum_x h(x).
\]
In particular, $\mu_d=\max_f\langle A_d\rangle_f$ is the largest mean
local field on the ordinary cube. We use $\|\cdot\|_{L^2}$ for the
norm from this inner product, $\|w\|_2^2=\sum_i w_i^2$ for the squared
Euclidean norm of the weights, and $I$ for the identity operator.

Each $\Gamma_i$ is symmetric and satisfies $\Gamma_i^2=I$ because its
edge sign is unchanged when bit $i$ is flipped. If $i<j$, flipping bit $i$
changes the sign in $\Gamma_j$, whereas flipping bit $j$ leaves the
sign in $\Gamma_i$ unchanged. The two paths across an $(i,j)$-square
therefore have opposite sign products, verifying full frustration and
giving $\Gamma_i\Gamma_j=-\Gamma_j\Gamma_i$. Flipping any bit also changes
parity, so $\Gamma_iP_d=-P_d\Gamma_i$. The mixed terms in
$\Gamma_w^2$ therefore cancel. For a Boolean function $f$ with local
field $h$, we have $h(x)^2=(\Gamma_w f)(x)^2$ at each vertex, so
\begin{equation}\label{eq:short-clifford-square}
 \Gamma_w^2=\|w\|_2^2I,\qquad
 \langle h,h\rangle=\langle f,\Gamma_w^2 f\rangle=\|w\|_2^2.
\end{equation}
The operator identity also appears in the Clifford-algebra formulation
of Mathews~\cite[Section~2.3]{Mathews} and in~\cite[Section~6.1]{LNSW}.

The two replacements below split weights among more
coordinates while preserving the sum of their squares. We refer to
them, in order, as the \emph{four-way refinement} and the
\emph{pair refinement}. A linear map $T$ between the
function spaces on two cubes \emph{preserves Boolean functions} if
$Tf$ takes its values in $\{-1,1\}$ whenever $f$ does. It is an
\emph{isometry} if $\langle Tf,Tg\rangle=\langle f,g\rangle$ for
all $f,g$, with the normalized inner product on each cube.

\begin{lemma}\label{lem:short-refinements}
Suppose that $w'$ is obtained from $w$ by one of the replacements
\[
 2b\longmapsto(b,b,b,b),\qquad
 (a,b)\longmapsto
 \left(\frac{a+b}{2},\frac{a+b}{2},
       \frac{a-b}{2},\frac{a-b}{2}\right),
\]
with $b>0$ in the first case and $a>b>0$ in the second. Write $d,d'$
for the lengths of $w,w'$. Then $\|w'\|_2^2=\|w\|_2^2$. There is
a linear isometry $T$ from real functions on $Q_d$ to real functions
on $Q_{d'}$ that preserves Boolean functions and satisfies
\[
 \Gamma_{w'} T=T\Gamma_w,\qquad P_{d'}T=TP_d.
\]
In particular, $\langle\Gamma_{w'}\rangle_{Tf}=\langle\Gamma_w\rangle_f$
for every Boolean function $f$ on $Q_d$.
\end{lemma}
The first identity allows us to apply the weighted adjacency operator
before or after refining. The second says that $T$ is
\emph{parity-preserving}. It ensures that refining one block of
coordinates respects the signs on edges in the directions that follow it.
\begin{proof}
The equality of the squared weight norms follows from
$(2b)^2=4b^2$ and
\[
 2\left(\frac{a+b}{2}\right)^2+
 2\left(\frac{a-b}{2}\right)^2=a^2+b^2.
\]
We first work on the one or two coordinates to be refined. At each
vertex of the four-bit cube we will select a vertex of the smaller
cube and copy the value of $f$ there, up to a sign.

Both constructions use a Boolean eigenfunction $u$ of $A_4$ with
eigenvalue $2$. To construct it, let $\eta$ be the indicator of
$\{000,111\}$ on $Q_3$, and put $g=A_3\eta$. The two vertices in
this set are not adjacent, and every vertex outside the set has
exactly one neighbor in it. Thus $g$ vanishes on the set and takes
values in $\{-1,1\}$ elsewhere. It follows that both $g+\eta$ and $g-\eta$
are Boolean. Since $A_3g=A_3^2\eta=3\eta$, direct multiplication gives
\[
 \begin{pmatrix}A_3&I\\ I&-A_3\end{pmatrix}
 \binom{g+\eta}{g-\eta}
 =\binom{2g+2\eta}{2g-2\eta}
 =2\binom{g+\eta}{g-\eta}.
\]
When the first bit indexes the two blocks, the displayed matrix is
$A_4$. The first coordinate contributes the off-diagonal identity
matrices, and the other three directions have opposite signs in the
two blocks. Hence $u=(g+\eta,g-\eta)$ is the required Boolean eigenfunction.

For the first replacement, let $\pi(x)$ be $0$ or $1$ according to
the parity of $|x|$, and define
\[
 (Tf)(x)=u(x)f(\pi(x))
\]
for a function $f$ of one bit. This map preserves Boolean functions,
since it multiplies a selected value of $f$ by a sign. On the one-bit
cube, the constant function $\mathbf1$ and $s\mapsto(-1)^s$ form a
basis of eigenfunctions of
$\Gamma_{(2b)}$ with eigenvalues $2b$ and $-2b$. Their images under
$T$ are $u$ and $P_4u$, with the same respective eigenvalues for $bA_4$,
because $A_4u=2u$ and $A_4P_4=-P_4A_4$. This proves
$bA_4T=T\Gamma_{(2b)}$. Moreover,
$(-1)^{\pi(x)}=(-1)^{|x|}$ gives
$P_4T=TP_1$ directly.

For the second replacement, let $\chi(s,t)=(-1)^s$ on the two-bit
source cube. Its weighted operator satisfies
\[
 \Gamma_{(a,b)}\mathbf1=a\mathbf1+b\chi,\qquad
 \Gamma_{(a,b)}\chi=b\mathbf1-a\chi.
\]
We seek a second Boolean function $v$ on four bits so that the refined
operator acts on $u,v$ by these same formulas. Group the four coordinate
operators in two pairs and set
\[
 L=\Gamma_0+\Gamma_1,\qquad R=\Gamma_2+\Gamma_3,
 \qquad v=(L-I)u.
\]
At each vertex $x$, the four quantities
$u(x)(\Gamma_i u)(x)$ are signs whose sum is
$u(x)(A_4u)(x)=2$. Exactly one of them is therefore negative.
In particular, the sum of the first two is either zero or two,
so
\[
 u(x)v(x)=u(x)(Lu)(x)-1\in\{-1,1\}.
\]
Since $u$ is Boolean, this proves that $v$ is Boolean as well.

The relations $L^2=R^2=2I$ and $RL=-LR$, together with $(L+R)u=2u$, yield
\[
 Lu=u+v,\qquad Ru=u-v,\qquad Lv=u-v,\qquad Rv=-u-v.
\]
For example, $Lv=L(L-I)u=2u-Lu=u-v$, and
$Rv=R(L-I)u=-LRu-Ru=-L(u-v)-(u-v)=-u-v$.
The refinement gives the left pair weight $(a+b)/2$ and the right pair
weight $(a-b)/2$, so $\Gamma_{w'}=\tfrac{a+b}{2}L+\tfrac{a-b}{2}R$. Hence
\[
 \Gamma_{w'}u=au+bv,\qquad \Gamma_{w'}v=bu-av.
\]
It remains to select the two source bits at each vertex $x$.
We want $T\mathbf1=u$ and $T\chi=v$, so the sign of the first bit
must be $u(x)v(x)$. The parity of $x$ then determines the second bit.
Thus we take
\[
 (-1)^{s(x)}=u(x)v(x),\qquad t(x)=s(x)\oplus\pi(x),
\]
and set
\[
 (Tf)(x)=u(x)f(s(x),t(x)).
\]
This map preserves Boolean functions. Since $(P_2\chi)(s,t)=(-1)^t$, the functions
$\mathbf1,\chi,P_2\chi,P_2\mathbf1$ form a basis on $Q_2$.
Their images under $T$, in this order, are
\[
 u,\quad v,\quad P_4v,\quad P_4u.
\]
Since $s(x)+t(x)\equiv\pi(x)\pmod2$, we also have $P_4T=TP_2$.
The matching formulas for $\Gamma_{w'}u$ and $\Gamma_{w'}v$ prove
$\Gamma_{w'} T=T\Gamma_{(a,b)}$ on
$\mathbf1,\chi$. The remaining two basis functions are their images under $P_2$.
Since both weighted operators anticommute with their parity operators
and $T$ preserves parity, the same identity holds on them as well.

Both maps select a source vertex $\phi(x)$ at each target vertex $x$.
The following argument applies to any map of the form
\[
 (Tf)(x)=\varepsilon(x)f(\phi(x)),\qquad \varepsilon(x)=\pm1,
\]
with $\Gamma_{w'}T=T\Gamma_w$ and positive source weights.
Let $T^*$ denote the adjoint of $T$ for the
normalized inner products. Grouping the target vertices according
to their value under $\phi$ shows that $T^*T$ is diagonal in the
vertex basis. Its entry at $y$ is
\[
 \rho(y)=2^{d-d'}\bigl|\phi^{-1}(y)\bigr|.
\]
The symmetry of both weighted operators and $\Gamma_{w'} T=T\Gamma_w$ give
\[
 \Gamma_w T^*T=T^*\Gamma_{w'} T=T^*T\Gamma_w.
\]
For adjacent source vertices $y,z$, the corresponding matrix entry
of this identity is
$\Gamma_w(y,z)\rho(z)=\rho(y)\Gamma_w(y,z)$. Since the source
weights are positive, $\rho(y)=\rho(z)$ on each edge. By connectivity
$\rho$ is constant, and its average is one because the sets
$\phi^{-1}(y)$ partition $Q_{d'}$. Hence $T^*T=I$
and $|\phi^{-1}(y)|=2^{d'-d}$ for every $y$.

To refine coordinates inside a larger cube, first move them to the
front. In the canonical signing a permutation of the bits also
requires a change of edge signs. If $\tau$ exchanges adjacent bits
$i,i+1$, define the signed swap
\[
 (Sf)(x)=(-1)^{x_ix_{i+1}}f(\tau x).
\]
It is a Boolean-preserving isometry with $S^2=I$. Substitution in
the definitions shows that
\[
 S\Gamma_iS=\Gamma_{i+1},\qquad
 S\Gamma_{i+1}S=\Gamma_i,
\]
while $S\Gamma_jS=\Gamma_j$ for $j\notin\{i,i+1\}$ and
$SP_d=P_dS$. Thus conjugating by $S$ exchanges the corresponding
weights. Repeating adjacent swaps places any selected coordinates
first.

Write $w=(\alpha,\beta)$, where $\alpha$ contains the weights of the
selected block and $\beta$ the remaining weights. If $\alpha$ has
length $k$, the weighted operator on the product of the two cubes is
\[
 \Gamma_w=\Gamma_\alpha\otimes I+P_k\otimes \Gamma_\beta.
\]
The factor $P_k$ occurs because the sign of every edge in the second
block depends on the parity of the first $k$ bits. If $\alpha'$ denotes
the refined weights and $k'$ their number, then $\Gamma_{\alpha'}T=T\Gamma_\alpha$ and
$P_{k'}T=TP_k$. Applying these identities to the two summands gives
$\Gamma_{(\alpha',\beta)}(T\otimes I)=(T\otimes I)\Gamma_w$. Likewise
$P_{d'}(T\otimes I)=(T\otimes I)P_d$, since parity on the full cube is
the product of the two block parities. Finally, $T\otimes I$ acts by
$T$ with the remaining bits fixed, so it is an isometry and preserves
Boolean functions. Preservation of the mean follows from
$\langle Tf,\Gamma_{w'}Tf\rangle=\langle Tf,T\Gamma_w f\rangle
=\langle f,\Gamma_w f\rangle$.
\end{proof}

When all the refined weights are one, $\Gamma_{w'}=A_n$. We call
a linear map $T$ from real functions on $Q_d$ to real functions on
$Q_n$ a \emph{lifting} of $w$ if it sends
Boolean functions to Boolean functions and satisfies $A_nT=T\Gamma_w$.

\begin{theorem}[Boolean lifting]\label{lem:short-unit-lifting}
Let $w=(w_0,\ldots,w_{d-1})$ have positive real entries. It admits a
lifting to $Q_n$ if and only if its entries are integers and
\[
 \|w\|_2^2=n,\qquad \gcd(w_0,\ldots,w_{d-1})=2^k
 \quad\text{for some integer } k\ge0.
\]
Every such lifting $T$ is an isometry. If the local field of a
Boolean function $f$ takes a given value at $m$ vertices of the
weighted cube $Q_d$, then the local field of $Tf$ takes that value
at $2^{n-d}m$ vertices of $Q_n$. In particular,
$\langle A_n\rangle_{Tf}=\langle\Gamma_w\rangle_f$.
One can choose $T$ so that $P_nT=TP_d$.
\end{theorem}
\begin{proof}
Suppose first that $w$ satisfies the stated conditions. Whenever
a weight is even, say $2b$, replace it by four copies of $b$.
If all weights are odd and two have distinct values $a>b$, apply
the second refinement of Lemma~\ref{lem:short-refinements} to $(a,b)$.
The new weights are again positive integers, and their squared norm
remains $n$.

Both refinements preserve the common odd divisors of the weights.
Indeed, an odd integer divides $2b$ if and only if it divides $b$.
If it divides $a$ and $b$, it divides their half-sum and half-difference,
since two is invertible modulo any odd integer. Conversely, the sum and
difference of the two new weights recover $a$ and $b$.
The greatest common odd divisor therefore remains one.

Each refinement increases the number of weights, which is at most $n$
since the weights are positive integers. The process therefore ends,
and it can end only when all weights have the same odd value. Their
greatest common odd divisor is one, so this value is one. The squared
norm then gives exactly $n$ weights. Composing the maps of
Lemma~\ref{lem:short-refinements} gives a lifting with $P_nT=TP_d$.

Conversely, let $T$ be any lifting of $w$. In any row of its matrix,
flipping one input sign changes the output by twice the corresponding
entry. Since both outputs belong to $\{-1,1\}$, every entry is in
$\{0,\pm1\}$. Choosing the input signs to agree with the nonzero
entries shows that their absolute values sum to one. Each row
therefore contains exactly one nonzero entry, and
\[
 (Tf)(x)=\varepsilon(x)f(\phi(x)),\qquad
 \varepsilon(x)\in\{-1,1\},
\]
for a map $\phi:V(Q_n)\to V(Q_d)$.

The entries of $A_nT$ are integers. In row $x$ of $T\Gamma_w$, the entries
at the distinct neighbors of $\phi(x)$ have magnitudes
$w_0,\ldots,w_{d-1}$. The intertwining identity consequently forces
all weights to be integers. Squaring it and using
\eqref{eq:short-clifford-square} gives
$nT=\|w\|_2^2T$, so $n=\|w\|_2^2$ because $T\ne0$.

To exclude an odd common divisor, observe that for every Boolean
function $s$ on $Q_n$ and every coordinate $j$ we have
\begin{equation}\label{eq:odd-divisor-obstruction}
 (\Gamma_js)^{\mathsf T}A_ns=2^n.
\end{equation}
Indeed, $\Gamma_jA_n=I+\sum_{i\ne j}\Gamma_j\Gamma_i$, and each
matrix $\Gamma_j\Gamma_i$ in the sum is skew-symmetric, so its
quadratic form vanishes. If an odd prime $p$ divided every weight,
the Boolean function $s=T\mathbf1$ would satisfy
$A_ns=T\Gamma_w\mathbf1\equiv0\pmod p$ pointwise. This contradicts
\eqref{eq:odd-divisor-obstruction}, whose left side would be divisible
by $p$. Thus the gcd is a power of two.

By the isometry argument in Lemma~\ref{lem:short-refinements},
$T^*T=I$ and each source vertex has $2^{n-d}$ preimages under $\phi$.
For a Boolean function $f$, the identity $A_nT=T\Gamma_w$ gives
\[
 (Tf)(x)(A_nTf)(x)
 =f(\phi(x))(\Gamma_w f)(\phi(x)).
\]
Thus the local field at each source vertex is repeated $2^{n-d}$ times
on $Q_n$.
\end{proof}

For positive integer weights, even realizing one local-field
distribution on an ordinary $Q_n$ requires $n=\|w\|_2^2$, by
\eqref{eq:short-clifford-square}. The gcd condition is necessary
as well. An odd prime dividing all the weights divides every local
field, so a Boolean function $s$ realizing the distribution would
have $A_ns\equiv0$ modulo that prime, contrary to
\eqref{eq:odd-divisor-obstruction}.

\section{Exact values from weighted cubes}\label{sec:exact}
\subsection{Attaining the signed-degree bound}
The equality condition in \eqref{eq:parity-bound} determines which
local fields we need to realize.

\begin{lemma}[Integer local-field bound]\label{lem:exact-local-field}
Let $f$ be a Boolean function on $Q_n$ and let
$h(x)=f(x)(A_nf)(x)$ be its local field. Then
\[
 h(x)\equiv n\pmod 2,\qquad
 \langle h,\mathbf1\rangle=\langle A_n\rangle_f,\qquad
 \langle h,h\rangle=n.
\]
Consequently, if $r$ is a positive integer of parity opposite to $n$,
then
\[
 \mu_n\le\frac{n+r^2-1}{2r}.
\]
A Boolean function attains this upper bound exactly when its local
field belongs to $\{r-1,r+1\}$ at every vertex. In that case the
proportion of vertices of field $r+1$ is $[n-(r-1)^2]/(4r)$
in each parity class of $Q_n$.
\end{lemma}
\begin{proof}
$h(x)$ is a sum of $n$ signs, so $h(x)\equiv n\pmod2$. Its mean
is $\langle A_n\rangle_f$ by definition, and its mean square is $n$
by \eqref{eq:short-clifford-square} with all weights equal to one.
If $r$ has parity opposite to $n$, then $h(x)-r$ is odd. Expanding
$(h(x)-r)^2\ge1$ and averaging over the vertices gives
\[
 1\le n-2r\langle A_n\rangle_f+r^2.
\]
Rearranging and maximizing over $f$ gives the bound on $\mu_n$, and
\eqref{eq:normalization} gives \eqref{eq:parity-bound}. Equality holds
exactly when $(h(x)-r)^2=1$ at every vertex. If a fraction $\alpha$
of the vertices then has field $r+1$, the second moment gives
\[
 n=(1-\alpha)(r-1)^2+\alpha(r+1)^2=(r-1)^2+4r\alpha.
\]
The two parity classes have the same size and the same sum of signed
degrees, since each edge contributes its sign once to each class.
With only two possible degrees, the proportions agree in both classes.
\end{proof}

\begin{proof}[Proof of \eqref{eq:exact-family} in Theorem~\ref{thm:exact-values}]
Put $n=r^2+1$. On the square with weights $(r,1)$, the constant
function $\mathbf1$ has local field $h(x)=r+(-1)^{x_0}$. The weights
are coprime and their squares sum to $n$, so
Theorem~\ref{lem:short-unit-lifting} lifts this function to $Q_n$
with the same two field values. Since $n$ and $r$ have opposite
parity, Lemma~\ref{lem:exact-local-field} gives equality in
\eqref{eq:parity-bound}, namely $F_n=(n-r)2^{n-2}$.
\end{proof}

For the six neighboring dimensions, start instead with a weighted
cube on which the local field is constant. If that constant is $r$,
changing one weight by $1$ changes each local field by $1$ or $-1$.
The change in the squared weight norm, and hence in the dimension
after lifting, depends on which weight we choose. The two cubes below
supply coordinates of weights $1$, $2$ and $3$.

\begin{lemma}\label{lem:exact-seeds}
For each of the weight vectors
\[
 (2,2,2,1,1,1,1),\qquad (3,1,1,1,1,1,1,1),
\]
the squared weights sum to $16$ and the weighted adjacency operator
has a Boolean eigenfunction with eigenvalue $4$.
\end{lemma}
\begin{proof}
Let $u$ be the Boolean function with $A_4u=2u$ constructed in
Lemma~\ref{lem:short-refinements}. Multiplying all four weights by
$2$ gives
\[
 \Gamma_{(2,2,2,2)}u=2A_4u=4u.
\]
Apply the first refinement to one coordinate of weight $2$,
replacing it by four coordinates of weight $1$. The intertwining
identity in Lemma~\ref{lem:short-refinements} shows that the resulting
Boolean function still has eigenvalue $4$. This gives the first vector.

For the second vector we use the same construction as for $A_4$ in
Lemma~\ref{lem:short-refinements}, with the length-seven Hamming
code in place of $\{000,111\}$. Let $H$ be the matrix whose columns
are the seven nonzero binary vectors of length three. Let $C$ consist of the vertices $x$
of $Q_7$ satisfying $Hx=0$ modulo two. Flipping bit $j$ adds the
$j$th column to $Hx$ modulo two. Thus a vertex in $C$ has no neighbor
in $C$. For a vertex outside $C$, exactly one column equals $Hx$,
so exactly one neighbor lies in $C$.

Let $\eta$ be the indicator of $C$ and put $g=A_7\eta$. On $C$,
$g$ vanishes. Outside $C$, its value is the sign of the unique edge
leading to $C$. Thus $g+\eta$ and $g-\eta$ are both Boolean.

Use a new first coordinate for the weight $3$, and write $w$ for
the second weight vector. Grouping vertices by this bit gives
\[
 \Gamma_w=\begin{pmatrix}A_7&3I\\3I&-A_7\end{pmatrix}.
\]
Since $A_7\eta=g$ and $A_7g=7\eta$, we obtain
\[
 \Gamma_w\binom{g+\eta}{g-\eta}
 =4\binom{g+\eta}{g-\eta}.
\]
The squared weights sum to $3\cdot2^2+4=16$ in the first vector and
$3^2+7=16$ in the second.
\end{proof}

\begin{proof}[Proof of \eqref{eq:additional-exact} in Theorem~\ref{thm:exact-values}]
Let $r$ be a positive integer divisible by $4$, and put $k=r/4$.
Fix $p\in\{\pm1,\pm2,\pm3\}$, and put $n=r^2+rp/2+1$.
This dimension is odd, so Lemma~\ref{lem:exact-local-field} gives
\[
 \mu_n\le\frac{n+r^2-1}{2r}=r+\frac p4.
\]
Choose a vector $w$ from Lemma~\ref{lem:exact-seeds} with
$w_j=|p|$, and let $f$ be its Boolean eigenfunction with eigenvalue
$4$. Multiplying the weights by $k$ makes the local field $4k=r$
and the squared weight norm $16k^2=r^2$.
Now change the weight $kw_j$ to $kw_j+\varepsilon$, where
$\varepsilon$ is the sign of $p$. Its square changes by $2kp+1$.
The resulting vector $w'$ therefore satisfies
\[
 \|w'\|_2^2=r^2+2kp+1=n,
 \qquad \Gamma_{w'}=k\Gamma_w+\varepsilon\Gamma_j.
\]
Since $\Gamma_w f=4f$ and $4k=r$, the local field of $f$ with the new weights is
\[
 f(x)(\Gamma_{w'}f)(x)
 =r+\varepsilon f(x)(\Gamma_jf)(x)
 \in\{r-1,r+1\}.
\]
Here $f(x)(\Gamma_jf)(x)$ is a sign, because both factors are signs.

Both vectors in Lemma~\ref{lem:exact-seeds} have at least two entries
equal to $1$, so $w'$ has an unchanged coordinate of weight $k$.
Its weights have gcd one because $k$ and $kw_j+\varepsilon$ are coprime.
The only zero weight occurs for $r=4$ and $p=-1$, giving
$n=15$. In this case restrict $f$ to the face $x_j=0$ and omit the
zero weight. The direction $j$ contributes nothing, and setting
$x_j=0$ leaves the signs of all remaining directions unchanged.
The restricted function still has local fields in $\{r-1,r+1\}$.

In every case Theorem~\ref{lem:short-unit-lifting} lifts the function
to $Q_n$ with local fields in $\{r-1,r+1\}$. Lemma~\ref{lem:exact-local-field}
therefore gives $\mu_n=r+p/4$, proving \eqref{eq:additional-exact}
by \eqref{eq:normalization}.
\end{proof}

\subsection{Boolean functions from optimal signings}\label{subsec:boolean-fields}
Let $r$ be a positive integer, put $n=r^2+1$, and fix an optimal
signing of $Q_n$. Its local field $h(x)$ is the signed degree at $x$.
Theorem~\ref{thm:exact-values} and Lemma~\ref{lem:exact-local-field}
give signed degree $r-1$ at half the vertices and $r+1$ at the other
half. Thus $g=h-r$ is a \emph{balanced} Boolean function. The choice
of lifting determines where these two values occur.

Fix a coordinate $i$. Its \emph{influence} on $g$ is the fraction
of vertices $x$ for which $g(x)\ne g(x\oplus e_i)$. Equivalently,
it is the fraction of edges in direction $i$ whose endpoints have
different signed degrees. Now fix a vertex $x$. The \emph{sensitivity}
of $g$ at $x$ is the number of coordinates $i$ for which the same
inequality holds, or equivalently the number of neighbors of $x$
whose signed degree differs from $h(x)$.

\begin{proposition}\label{prop:field-directions}
Let $r\ge3$ be odd and put $n=r^2+1$. There is an optimal signing
of $Q_n$ with local field $h$ such that every coordinate has influence
$1/2$ on $g=h-r$ and $g$ has sensitivity $n/2$ at every vertex of
odd Hamming weight.
\end{proposition}
\begin{proof}
Put $a=(r+1)/2$, $b=(r-1)/2$ and $m=a^2+b^2=n/2$. Use the weighted
square from the proof of \eqref{eq:exact-family}. Every lift gives an
optimal signing, and $g$ is the sign $(-1)^{x_0}$ of the first
source bit selected by the lift.
First refine $(r,1)$ to $(a,a,b,b)$ and group the coordinates into
two pairs of weights $(a,b)$. Since $a$ and $b$ are consecutive,
Theorem~\ref{lem:short-unit-lifting} gives a parity-preserving lift
of each pair to $Q_m$. Use the same lift for both pairs, and write
$X,Y\in Q_m$ for the two blocks of input bits. Let
$\phi\colon Q_m\to Q_2$ be the map selecting the two source bits
in this lift. Each of the four possible pairs is selected equally
often, and
$|\phi(X)|\equiv |X|\pmod 2$.

In the refinement the four bits occur in the order $(s,s',t,t')$,
of weights $(a,a,b,b)$. The two blocks group them as
$(s,t)=\phi(X)$ and $(s',t')=\phi(Y)$. Reordering coordinates by the
signed swaps of Lemma~\ref{lem:short-refinements} only permutes the
local fields. In the table the first block selects the row, the second
selects the column, and the entry is $g(X,Y)$.
\[
 \begin{array}{c|rrrr}
 (s,t)\backslash(s',t')&00&01&10&11\\ \hline
 00&-1&-1& 1&-1\\
 01& 1& 1& 1&-1\\
 10&-1& 1& 1& 1\\
 11&-1& 1&-1&-1
 \end{array}
\]
Indeed, the first source bit has sign $uv$ in the pair refinement of
Lemma~\ref{lem:short-refinements}. Substituting its function $u$ in
$v=(\Gamma_0+\Gamma_1-I)u$ gives $uv=-(-1)^{s'+t}$ when $t=t'$
and $-(-1)^{s+t}$ when $t\ne t'$. The entry is $-1$ precisely when
$(s,t,s',t')$ is nondecreasing or nonincreasing. Thus $g(X,Y)$ is
obtained by applying the four-bit function of Ambainis~\cite{Ambainis}
to the selected bits. The same four-bit function occurs in the vertex
coloring of Laplante et al.~\cite[Section~5.3]{LNSW}.

Fix a bit of $X$ and a value of $X$. Flipping that bit changes the
parity of $X$, so it selects a different row. The two rows disagree
in exactly two columns, and as $Y$ varies each column is selected
equally often. Thus this bit has influence $1/2$. Interchanging rows
and columns gives the same conclusion for every bit of $Y$.

To count sensitivity, observe that at an even-weight vertex the row
and column have the same parity.
The table entry is then $-1$ for even rows and $1$ for odd rows.
From an odd-weight vertex, flipping a bit of $X$ changes the row
parity and flipping a bit of $Y$ preserves it. All these neighbors
are even, with $m$ of each row parity. Hence $m=n/2$ neighbors have
each sign, and exactly $n/2$ have sign different from the vertex itself.
\end{proof}

\begin{corollary}\label{cor:field-every-vertex}
Let $r=2^k\pm1$ for an integer $k\ge2$, and put $n=r^2+1$. There
is an optimal signing of $Q_n$ with local field $h$ such that every
coordinate has influence $1/2$ on $g=h-r$ and $g$ has sensitivity
$n/2$ at every vertex. Thus every vertex has $n/2$ neighbors of
signed degree $r-1$ and $n/2$ of signed degree $r+1$.
\end{corollary}
\begin{proof}
We choose the lift of $(r,1)$ in Proposition~\ref{prop:field-directions}
so that at every vertex exactly half the incident edges change the
first selected bit. We induct on $r$, starting with the identity
on the two-bit cube for $(1,1)$.

Write $r=2^k+\varepsilon$ with $k\ge2$ and $\varepsilon\in\{-1,1\}$,
using $3=4-1$. The two weights $a,b$ in the proposition are
$u=2^{k-1}$ and $v=u+\varepsilon$, in some order. Split $u$ into
$u^2$ coordinates of weight one by repeated four-way refinements.
Pair one of these coordinates with $v$. Since $v=2^{k-1}+\varepsilon<r$,
we may use the inductive lift of $(v,1)$ (with $v=1$ in the base
case). Exactly $(v^2+1)/2$ flips change the original bit of weight
$v$, all within this lift. The full block has $m=u^2+v^2=n/2$ bits.
Each flip changes the parity of the two selected bits, so it changes exactly
one selected bit. The remaining $m-(v^2+1)/2$ flips change the bit
of weight $u$. Both counts are independent of the vertex.

Use this lift in both blocks of the proposition. If $d$ flips of
$X$ change $s$, the other $m-d$ change $t$. The
same counts hold for $s',t'$ in $Y$. In the table above, flipping
$s$ and $s'$ together reverses the sign, as does flipping $t$ and
$t'$ together. The two individual flips in either pair therefore
lead to opposite entries, and exactly one changes the sign. At every
vertex of $Q_n$, the number of flips changing the table entry is
therefore $d+(m-d)=m$.
The table entry is the sign of the first bit on $(r,1)$, so this
proves the induction claim and gives sensitivity $m$ everywhere.
Optimality and the coordinate influences follow from
Proposition~\ref{prop:field-directions} with this choice of lift.
\end{proof}

\section{The uniform bound}\label{sec:asymptotic}
To pass from the weighted construction to an arbitrary dimension $n$,
we use the fact that for every $d\ge1$,
\begin{equation}\label{eq:monotone}
 \mu_{d+1}\ge\mu_d.
\end{equation}
Indeed, take two copies of a fully frustrated signing attaining $\mu_d$.
Join each vertex $x$ in one copy to its counterpart in the other by an
edge of sign $(-1)^{|x|}$. Each new square has two opposite joining signs
and two equal signs from the copies. It is therefore negative. Half the
joining edges are positive, so the mean local field is still $\mu_d$.

For $d\le n$, it follows that
\[
 \sqrt n-\mu_n\le
 \frac{n-d}{\sqrt n+\sqrt d}+(\sqrt d-\mu_d).
\]
Thus it suffices to keep $n-d$ bounded and to have
$\sqrt d-\mu_d=O(d^{-1/2})$. When our starting dimension $r^2+1$ is
too far below $n$, we change one weight to reach a dimension within
$15$ of $n$. For this we first need a coordinate of a prescribed weight.

\begin{lemma}\label{lem:one-weight}
Let $q$ and $a$ be integers with $1\le a\le q/2$. Starting with two
coordinates of weight $q$ and $q$ coordinates of weight $1$, the
refinements of Lemma~\ref{lem:short-refinements} can produce a
coordinate of weight $a$.
\end{lemma}
\begin{proof}
We induct on $q$. If $a=1$, the desired weight is already present.
Otherwise, keep one copy of $q$ and split the other into two copies
each of $u=\lfloor q/2\rfloor$ and $v=\lceil q/2\rceil$. For even
$q$ this is the four-way refinement. For odd $q$ refine the pair
$(q,1)$. Set aside two copies of $v$.

We seek two weights whose half-sum or half-difference is $a$. Choose
$c$ closest to $2a$ among $0,u,q$, and put $b=|2a-c|$. Since
$2a\in[0,q]$ and the two gaps in the list have lengths $u,v$, we
have $b\le v/2$. If $b>0$,
the induction hypothesis gives a coordinate of weight $b$ from the
two copies of $v$ and at most $v$ further coordinates of weight $1$.
These are available since $1+v\le q$. If $c>0$, we also have a
separate coordinate of weight $c$.

If $b,c$ are positive and distinct, they have the same parity and
$a$ is either $(b+c)/2$ or $|b-c|/2$. Both weights occur in the
pair refinement. If $b=c>0$, then $c=a$.
If one of $b,c$ is zero, the other is $2a$ and can be split into
four copies of $a$.
\end{proof}

\begin{proof}[Proof of Theorem~\ref{thm:main}]
It suffices to consider $n\ge5$. Let $r=2q$ be the largest even
integer with $r^2+1\le n$, and put $\Delta=n-r^2-1$. Then
$n\le(r+2)^2$, so $\Delta\le8q+3$. If $\Delta<80$,
\eqref{eq:monotone} and \eqref{eq:exact-family} give
\[
 \sqrt n-\mu_n\le\sqrt n-r
 =\frac{\Delta+1}{\sqrt n+r}\le\frac{40}{r}.
\]

Suppose now that $\Delta\ge80$. Increasing a weight $a$ by eight
adds $16a+64$ to the squared weight norm. These increments are
spaced by $16$, so we choose the largest one at most $\Delta$. Put
\[
 a=\left\lfloor\frac{\Delta-64}{16}\right\rfloor,\qquad
 d=r^2+1+16a+64.
\]
Since $80\le\Delta\le8q+3$, we have $1\le a\le q/2$.
The choice of $d$ gives $d\le n\le d+15$. Start with the
constant function on the cube with weights $(r,1)$, whose mean local
field is $r$. Split $r$ into four copies of $q$. The block $(q,q,1)$
has gcd one, so the refinements in the proof of
Theorem~\ref{lem:short-unit-lifting} turn it into $2q^2+1$ coordinates
of weight $1$. Set one aside to keep the gcd equal to one after
the weight change. Use
Lemma~\ref{lem:one-weight} on the other two copies of $q$ and $q$ of
the remaining coordinates of weight $1$ to obtain the weight $a$.

Let $w$ and $f$ be the resulting weights and Boolean function, with
$j$ the coordinate of weight $a$ obtained from the lemma. The refinements
preserve $\|w\|_2^2=r^2+1$ and $\langle\Gamma_w\rangle_f=r$.
Increasing $w_j$ by eight adds $8\langle\Gamma_j\rangle_f$ to the
mean local field. To estimate this contribution, observe that
\eqref{eq:short-clifford-square} gives
\[
 \|(\Gamma_w-rI)f\|_{L^2}^2=(r^2+1)-2r^2+r^2=1.
\]
For $i\ne j$, $\Gamma_j\Gamma_i$ is skew-symmetric, so
$\langle\Gamma_j f,\Gamma_w f\rangle=w_j=a$. As
$\|\Gamma_j f\|_{L^2}=1$, Cauchy--Schwarz gives
\[
 \big|a-r\langle\Gamma_j\rangle_f\big|
 =\big|\langle\Gamma_j f,(\Gamma_w-rI)f\rangle\big|\le1.
\]
Thus the increase in the mean is at least $8(a-1)/r$. The new squared weight norm
is $d$, and the coordinate of weight $1$ kept above ensures that
the new weights have gcd one. Theorem~\ref{lem:short-unit-lifting}
and \eqref{eq:monotone} now give
\[
 \mu_n\ge\mu_d\ge r+\frac{8(a-1)}{r}.
\]
On the other hand, $n-r^2\le16a+80$, so
$\sqrt n-r\le(n-r^2)/(2r)\le(8a+40)/r$. Subtracting the two
bounds gives $\sqrt n-\mu_n\le48/r$. Since $r\ge\sqrt n-2$,
the theorem follows from \eqref{eq:normalization}, after enlarging
the absolute constant for $n<5$.
\end{proof}

\section{Integral eigenvectors and dimension nine}\label{sec:q9}
Let $n=r^2$ with $r$ a positive integer. A constant signed degree
of $r$ now satisfies both the parity and the second-moment conditions
of Lemma~\ref{lem:exact-local-field}. To improve the spectral bound,
we measure how far a switching is from having this constant degree.
Suppose a Boolean function $s$ gives $(n-r)2^{n-2}+\delta$ negative
edges, and let $q=(A_ns-rs)/2$ be its residual.
At a vertex $x$ with local field $h(x)$ we
have $q_x=s(x)(h(x)-r)/2$. Since $h(x)\equiv n\equiv r\pmod2$,
$q$ has integer entries, and $A_n^2=r^2I$ gives $A_nq=-rq$.
Counting each edge twice with its switched sign gives
$\sum_x h(x)=r2^n-4\delta$, while Lemma~\ref{lem:exact-local-field}
gives $\sum_x h(x)^2=2^n r^2$. Hence in the unnormalized Euclidean norm
\begin{equation}\label{eq:square-residual}
 \|q\|_2^2=\tfrac14\sum_x(h(x)-r)^2
 =\tfrac14\Bigl(2^{n+1}r^2-2r\sum_x h(x)\Bigr)=2r\delta.
\end{equation}
For $Q_9$, ruling out $768$ through $771$ negative edges means
excluding squared norms $0,6,12,18$, respectively. If $q=0$, then
$A_9s=3s$, which would force $3\mid2^9$ by
\eqref{eq:odd-divisor-obstruction}. We shall show that the least
positive squared norm of an integral $-3$-eigenvector of $A_9$ is $18$.
To exclude $771$ edges, we must also identify the vectors of that norm
and rule them out as residuals of Boolean switchings.

We need two estimates for a vector supported on a small set of
vertices. The first bounds the signed adjacency operator restricted
to that set. Write
$\|B\|_{\rm op}=\max_{\|v\|_2=1}\|Bv\|_2$ for the Euclidean operator norm.

\begin{lemma}\label{lem:principal}
Let $S$ be a nonempty set of vertices of $Q_d$, and let $A_d[S]$ be
the principal submatrix on $S$. Then $\|A_d[S]\|_{\rm op}\le\sqrt{|S|-1}$.
\end{lemma}
\begin{proof}
Write $B=A_d[S]$. The entry $(B^2)_{xy}$ is the sum of the signed products along
length-two walks from $x$ to $y$ whose intermediate vertex belongs to $S$.
For $x\ne y$, such walks exist only if $x,y$ have distance two.
There are exactly two possible intermediate vertices. If both belong
to $S$, the two products cancel because their square is negative.
If just one belongs to $S$, the entry has magnitude one.

Fix a vertex $x\in S$. The diagonal entry $(B^2)_{xx}$ counts its
neighbors in $S$, all of opposite parity to $x$. Nonzero off-diagonal
entries in row $x$ occur only at other vertices of the same parity
as $x$, since two flips preserve parity. Each has magnitude one.
Thus the absolute row sum is at most $|S|-1$.
The spectral radius of a matrix is at most its maximum absolute row
sum. Since $B$ is symmetric, $\|B\|_{\rm op}^2$ is the largest
eigenvalue of $B^2$. Hence $\|B\|_{\rm op}^2\le |S|-1$.
\end{proof}

The second estimate bounds the number of negative edges inside the
support. They form a bipartite graph, and they contain no $4$-cycle
because four negative edges around a square would have positive sign
product.

\begin{lemma}\label{lem:bipartite-c4}
A bipartite $C_4$-free graph on $t\ge2$ vertices has at most
$t(1+\sqrt{2t-3})/4$ edges.
\end{lemma}
\begin{proof}
Let the two parts have sizes $a$ and $b$, with $a+b=t$, and let $e$
be the number of edges. Any two vertices in the first part have at most
one common neighbor. Summing over the $b$ vertices $v$ in the second
part gives $\sum_v\binom{\deg(v)}2\le\binom a2$.
Since $\sum_v\deg(v)=e$, Cauchy--Schwarz gives
\[
 e^2\le b\sum_v\deg(v)^2\le b\bigl(a(a-1)+e\bigr).
\]
Interchanging the parts gives $e^2\le a\bigl(b(b-1)+e\bigr)$.
Adding and using $ab\le t^2/4$ yields
\[
 2e^2\le ab(t-2)+te\le t^2(t-2)/4+te.
\]
Solving this quadratic inequality for $e$ proves the assertion.
\end{proof}

The next lemma identifies the shortest nonzero integral vectors in
the negative eigenspace. The vectors in its equality statement are
supported on a vertex and its neighbors. Such eigenvectors appear
in~\cite[Section~6.3]{LNSW}. Write $\mathbf e_z$ for the vector with
entry one at vertex $z$ and zero elsewhere.

\begin{lemma}\label{lem:integral-minimum}
Let $n=r^2$ with $r\ge3$ an integer. Every nonzero integral vector $q$ with
$A_nq=-rq$ satisfies $\|q\|_2^2\ge2r^2$. Equality holds precisely when
$q=\pm(A_n-rI)\mathbf e_z$ for some vertex $z$.
\end{lemma}
\begin{proof}
We examine the half of the cube carrying less of the squared norm of
$q$. The eigenvector equation determines the other half, so it is
enough to control the support and entries on this one. Ordering vertices
by their first coordinate gives
\[
 A_n=\begin{pmatrix}A_{n-1}&I\\ I&-A_{n-1}\end{pmatrix}.
\]
Write $q=(u,v)$, interchanging the two halves if necessary so that
$k=\|u\|_2^2\le\|v\|_2^2$. After this interchange the block matrix
has the form $\left(\begin{smallmatrix}C&I\\I&-C\end{smallmatrix}\right)$,
where $C=A_{n-1}$ or $-A_{n-1}$, and $C^2=(r^2-1)I$.
The first block of the eigenvector equation reads $Cu+v=-ru$, so
$v=-ru-Cu$. Neither half is zero, since either zero half would force
the other to vanish. Expanding the squared norm gives
\[
 \begin{aligned}
 \|q\|_2^2
 &=k+\|ru+Cu\|_2^2\\
 &=k+r^2k+2r\,u^{\mathsf T}Cu+(r^2-1)k
 =2r\bigl(rk+u^{\mathsf T}Cu\bigr).
 \end{aligned}
\]
Set $a=rk+u^{\mathsf T}Cu$. This is a positive integer, since $u$ and $C$
have integer entries and $\|q\|_2^2=2ra>0$. As $u$ was the half of
smaller norm, $2k\le\|q\|_2^2$. Thus
\begin{equation}\label{eq:integral-norm}
 \|q\|_2^2=2ra,\qquad k\le ra.
\end{equation}
Let $S$ be the support of $u$, and put $t=|S|$. Every nonzero integral
entry contributes at least one to $k$, so $t\le k$.
Only the principal submatrix $C[S]$ contributes to $u^{\mathsf T}Cu$.
Lemma~\ref{lem:principal}, which also applies to $-A_{n-1}$, gives
\[
 |u^{\mathsf T}Cu|\le\|C[S]\|_{\rm op}\|u\|_2^2
 \le k\sqrt{t-1}\le k\sqrt{k-1}.
\]
It follows that
\begin{equation}\label{eq:small-support}
 a\ge k\bigl(r-\sqrt{k-1}\bigr).
\end{equation}

Suppose $\|q\|_2^2\le 2r^2$, or equivalently $a\le r$. It suffices
to show that $q$ then has the asserted form and attains equality.
Equation~\eqref{eq:integral-norm} gives $1\le k\le ra\le r^2=n$.
For $2\le k\le n-2$, we have
\[
 k\bigl(r-\sqrt{k-1}\bigr)>r
 \quad\Longleftrightarrow\quad k^2-nk+n<0.
\]
Indeed, rearranging and dividing by $\sqrt{k-1}$ reduces the first
inequality to $r\sqrt{k-1}>k$, which can be squared.
The convex quadratic $k^2-nk+n$ has value $4-n<0$ at both endpoints
of this interval, and is therefore negative throughout it.
Equation~\eqref{eq:small-support} would then imply $a>r$.
Thus $k\in\{1,n-1,n\}$.

If $k=1$, the integral vector $u$ has a single nonzero entry.
Thus $u=\varepsilon\mathbf e_y$ for some vertex $y$ of $Q_{n-1}$
and some sign $\varepsilon$.
Since $C$ has zero diagonal, $u^{\mathsf T}Cu=0$ and $a=r$.
Let $z$ be the vertex in the $v$-half with the same last $n-1$ bits
as $y$. Then
$\mathbf e_z=(0,\mathbf e_y)$. The formula for $v$ gives
\[
 q=\varepsilon\binom{\mathbf e_y}{-r\mathbf e_y-C\mathbf e_y}
   =\varepsilon(A_n-rI)\mathbf e_z.
\]
It remains to exclude $k=n-1,n$. In these two cases the bound
$a\le r$ forces so many negative edges on the support of $u$ that
Lemma~\ref{lem:bipartite-c4} will give a contradiction. Indeed,
\begin{equation}\label{eq:large-negative-form}
 u^{\mathsf T}Cu=a-rk\le-r(k-1).
\end{equation}
If an entry of $u$ has magnitude at least two, then
$k\ge4+(t-1)$, so $t\le k-3$. Lemma~\ref{lem:principal} would give
\[
 |u^{\mathsf T}Cu|\le k\sqrt{k-4}<r(k-1),
\]
contrary to \eqref{eq:large-negative-form}. The strict inequality
follows by squaring, since
$k^2(k-4)<k(k-1)^2\le r^2(k-1)^2$.
Thus the support has exactly $k$ vertices, and every nonzero entry
of $u$ is $1$ or $-1$.

Switch the induced signed graph on this support by the signs $u_x$.
Each negative edge then contributes $-2$ to $u^{\mathsf T}Cu$, and
each positive edge contributes $2$. Thus \eqref{eq:large-negative-form}
requires at least $r(k-1)/2$ negative edges. Switching preserves the negative
sign product on each square, so the negative edges on the support
still form a bipartite $C_4$-free graph.
Comparing with Lemma~\ref{lem:bipartite-c4} and multiplying by $4/k$
gives $2r(1-1/k)\le1+\sqrt{2k-3}$.
But $r^2-1\le k\le r^2$ and $r\ge3$ give $k\ge2r$ and
\[
 2r(1-1/k)\ge2r-1
 >1+\sqrt{2r^2-3}\ge1+\sqrt{2k-3},
\]
a contradiction. The strict middle inequality follows from
$(2r-2)^2-(2r^2-3)=2(r-2)^2-1>0$.

Thus $k=1$ is the only possibility when $\|q\|_2^2\le 2r^2$.
In this case we found $a=r$ and $q=\pm(A_n-rI)\mathbf e_z$.
Conversely, for every $z$ we have
\[
 A_n(A_n-rI)\mathbf e_z=(r^2I-rA_n)\mathbf e_z
 =-r(A_n-rI)\mathbf e_z.
\]
The vector $(A_n-rI)\mathbf e_z$ has entry $-r$ at $z$, an entry of magnitude
one at each of its $n$ neighbors, and zero elsewhere. Its squared
norm is $r^2+n=2r^2$, as required.
\end{proof}

It remains to decide whether a shortest integral eigenvector can be
the residual of a switching. The case $q=0$ includes the classical
construction at powers of four~\cite{BHN,LNSW}.

\begin{proposition}\label{prop:square-gap}
Equality in the spectral bound $F_n\ge(n-\sqrt n)2^{n-2}$ holds
if and only if $n$ is a power of four. If $r\ge3$ is an integer
that is not a power of two, then
\begin{equation}\label{eq:square-gap}
 F_{r^2}\ge(r^2-r)2^{r^2-2}+
 \begin{cases}
 r+1,&r\text{ odd},\\
 r+2,&r\text{ even}.
 \end{cases}
\end{equation}
\end{proposition}
\begin{proof}
Since $A_n$ is symmetric and $A_n^2=nI$, equality in its Rayleigh
bound forces a Boolean vector $s$ with $A_ns=\sqrt n\,s$.
Every entry of $A_ns$ is integral and every entry of $s$ is $\pm1$,
so $r=\sqrt n$ must be an integer.
Equation~\eqref{eq:odd-divisor-obstruction} now gives
$2^n=(\Gamma_i s)^{\mathsf T}A_ns=r(\Gamma_i s)^{\mathsf T}s$
for any coordinate $i$. Since $(\Gamma_i s)^{\mathsf T}s$ is an
integer, $r$ divides $2^n$ and is a power of two. Conversely, the constant
function on the one-edge cube of weight $r$ has local field $r$
at both vertices. If $r$ is a power of two,
Theorem~\ref{lem:short-unit-lifting} lifts this function to $Q_{r^2}$
with the same constant field, attaining the spectral bound.

Suppose now that $r\ge3$ is not a power of two, and choose $s$ attaining
$F_n$ for $n=r^2$. Then $q\ne0$ and $\delta=F_n-(n-r)2^{n-2}$.
By \eqref{eq:square-residual} and Lemma~\ref{lem:integral-minimum},
$2r\delta=\|q\|_2^2\ge2r^2$, so $\delta\ge r$.

Suppose $\delta=r$. The equality case of the lemma gives
$q=\varepsilon(A_n-rI)\mathbf e_z$ for some vertex $z$ and sign $\varepsilon$.
Thus $(A_n-rI)s=2\varepsilon(A_n-rI)\mathbf e_z$. Subtracting
$2\varepsilon\mathbf e_z$ from $s$ gives a vector $s'$ satisfying
$A_ns'=rs'$. Only its entry at $z$ could fail to be a sign.

The squared entries of $s'$ sum to the same value on the two parity
classes. Indeed, $A_nP_n=-P_nA_n$ puts $P_ns'$ in the $-r$ eigenspace,
orthogonal to $s'$, and $(s')^{\mathsf T}P_ns'$ is the difference
of the two sums. Each parity class has $2^{n-1}$ vertices, and all
entries of $s'$ except possibly the one at $z$ have square one.
Equality of the sums therefore forces $(s'_z)^2=1$ as well.
Thus $s'$ is Boolean, and its eigenvector equation
contradicts the spectral nonattainment already proved. Hence $\delta>r$.

Since $\delta$ is integral, we obtain $\delta\ge r+1$.
If $r$ is even, then $n=r^2\ge16$ is even, so switching at a single vertex
changes the number of negative edges by an even number.
In the canonical signing, the first direction is positive and each
other direction has $2^{n-2}$ negative edges. Its total
$(n-1)2^{n-2}$ is even, so every switching has an even number of
negative edges, as does the spectral baseline
$(n-r)2^{n-2}$. It follows that $\delta$ is even and
$\delta\ge r+2$ in this case.
\end{proof}

\begin{proof}[Proof of Theorem~\ref{thm:nine}]
Putting $r=3$ in Proposition~\ref{prop:square-gap} gives
$F_9\ge(9-3)2^7+4=772$.
The explicit switching in~\cite{Q9witness} has $772$ negative edges,
giving equality.
\end{proof}

\section{Concluding remarks}\label{sec:conclusion}
\noindent\textbf{Square-free color classes.}\enspace
Let $B_n$ be the largest size of a color class in an edge two-coloring
of $Q_n$ with no monochromatic square. The construction in
Corollary~\ref{cor:extremal} gives $B_n\ge n2^{n-1}-F_n$.
For the reverse bound in the dimensions
of Theorem~\ref{thm:exact-values}, we must also allow squares with two
edges of each color.

\begin{proposition}\label{prop:square-count}
Let $T\subseteq E(Q_n)$ meet every square, and let $r$ be a positive
integer of parity opposite to $n$. Write $M$ for the number of squares
whose four edges belong to $T$. A square is \emph{alternating} if
it meets $T$ in exactly two opposite edges. If there are $A$ such
squares, then
\[
 |T|\ge\left(n-\frac{n+r^2-1}{2r}\right)2^{n-2}+\frac{A-M}{r}.
\]
With the edges of $T$ signed negative and the others positive,
equality holds if and only if every signed degree belongs to
$\{r-1,r+1\}$.
\end{proposition}
\begin{proof}
Write $h(x)=n-2\deg_T(x)$, so $\sum_x h(x)=n2^n-4|T|$.
In $\sum_x h(x)^2$ the squared signs contribute $n2^n$.
Each pair of distinct incident edges lies on a unique square.
For a square with signs $a,b,c,d$ in cyclic order, the cross terms
sum to $2(ab+bc+cd+da)=2(a+c)(b+d)$. Since $T$ meets every square,
this is $8$ for a square contained in $T$, $-8$ for an alternating
square and zero otherwise. Thus
\[
 \sum_x h(x)^2=n2^n+8(M-A).
\]
Since each $h(x)-r$ is odd, $\sum_x(h(x)-r)^2\ge2^n$, with equality
exactly when $|h(x)-r|=1$ at every vertex. Substituting the two sums
gives the claimed bound.
\end{proof}

For the smaller class $T$ of a coloring attaining $B_n$, we have
$M=0$ and $|T|=n2^{n-1}-B_n$. In the dimensions of
Theorem~\ref{thm:exact-values} the proposition gives $|T|\ge F_n+A/r$.
Thus $B_n=n2^{n-1}-F_n$, and every maximizing coloring in these
dimensions has no alternating square. For general $n$,
Corollary~\ref{cor:extremal} and Bialostocki's upper bound~\cite{Bialostocki}
give
$B_n=(n+\sqrt n+O(n^{-1/2}))2^{n-2}$.

\smallskip\noindent\textbf{The uniform constant.}\enspace
Since $\mu_2=1$, any constant in Theorem~\ref{thm:main} satisfies
$C\ge2-\sqrt2$. Refining the estimates for the weighted construction
and checking the smaller dimensions gives an admissible $C<1$.
Whether $C=2-\sqrt2$ suffices remains open.

\smallskip\noindent\textbf{Dimension 25.}\enspace
The next square dimension in which the spectral bound is strict is
$25$. The exact value remains open. Beyond the bound of
Proposition~\ref{prop:square-gap}, we record without proof
\[
 167\,772\,180\ \le F_{25}\le\ 168\,230\,912.
\]
The lower endpoint uses restrictions on the first three distance
layers of an integral eigenvector of $A_{25}$. The upper endpoint
comes from a switching constructed with three copies of the binary
Hamming code of length seven.

\section*{Acknowledgments}
The author thanks Reza Naserasr for valuable suggestions that
improved the presentation and for raising the
question of which Boolean functions arise from optimal signings.
The author is partially supported by the National Natural Science
Foundation of China (grant No.~12371359).

\section*{Declaration of generative AI and AI-assisted technologies}
During the preparation of this work, the author used OpenAI Codex to
explore proof strategies, check mathematical arguments, run computational
experiments, search the literature, and assist with drafting and revising
the manuscript. The author reviewed and revised the resulting material
and takes full responsibility for the content of the article.

\end{document}